\documentclass[12pt]{amsart}
\usepackage{amsfonts}
\usepackage{amssymb, amstext, amscd, amsmath}
\usepackage{amsthm}
\usepackage{times}
\usepackage{txfonts}
\usepackage{indentfirst}
\usepackage[all]{xy}
\usepackage{subfigure}
\usepackage{graphicx}
\newtheorem{thm}{Theorem}[section]

\newtheorem{lem}[thm]{Lemma}
\newtheorem{conj}[thm]{Conjecture}

\theoremstyle{definition}
\newtheorem{rem}[thm]{Remark}

\newtheorem{dfn}[thm]{Definition}

\numberwithin{equation}{section}
\numberwithin{equation}{section}
\def\square{\vbox{
      \hrule height 0.4pt
      \hbox{\vrule width 0.4pt height 5.5pt \kern 5.5pt \vrule width 0.4pt}
      \hrule height 0.4pt}}

\def\ch\mathrm{c h}

\begin{document}

\title[The Slope Conjecture for 3-String Montesinos Knots]
{The Slope Conjecture for 3-String Montesinos Knots}

\author{Xudong Leng$^{*}$}
\address{(Xudong Leng) School of Mathematical Sciences, Dalian University of
Technology, Dalian 116024, P. R. China} \email{xudleng@163.com}
\thanks{$^{*}$ Corresponding author}
\author{Zhiqing Yang$^{\dag}$}
\address{(Zhiqing Yang) School of Mathematical Sciences, Dalian University of
Technology, Dalian 116024, P. R. China} \email{yangzhq@dlut.edu.cn}
\thanks{$^{\dag}$ Supported by the NFSC (No. 11271058)}
\author{Xinmin Liu$^{\ddag}$}
\address{(Xinmin Liu) School of Mathematical Sciences, Dalian University of
Technology, Dalian 116024, P. R. China }
\email{ximinliu@dlut.edu.cn}
\thanks{$^{\ddag}$ Supported by the NFSC (No. 11431009) }



\begin{abstract}
The (Strong) Slope Conjecture relates the degree of the colored Jones polynomial of a knot to certain essential surfaces in the knot complement. We verify the Slope Conjecture and the Strong Slope Conjecture for 3-string Montesinos knots satisfying certain conditions.
\end{abstract}
\maketitle


\section{Introduction}
The colored Jones polynomial is a generalization of the Jones polynomial, a celebrated knot invariant originally discovered by V. Jones via von Neumann algebras ~\cite{J}. To find an intrinsic interpretation of the Jones polynomial, E. Witten introduced the Chern-Simons quantum field theory ~\cite{W} which led to invariants of 3-manifolds as well as the colored Jones polynomial. Then N. Reshetikhin and V. Turaev constructed a mathematically rigorous mechanism by quantum groups ~\cite{RT} to produce these invariants. Later, G. Masbaum and P. Vogel defined the colored Jones polynomial through skein theory ~\cite{MV} using the Tiemperley-Lieb algebra.

 Compared with the Jones polynomial, the colored Jones polynomial reveals much stronger connections between quantum algebra and 3-dimensional topology, for example, the Volume Conjecture, which relates the asymptotic behavior of the colored Jones polynomial of a knot to the hyperbolic volume of its complement. Another connection proposed by S. Garoufalidis ~\cite{Gar11b} named Slope Conjecture, predicts that the growth of maximal degree of the colored Jones polynomial of a knot determines some boundary slopes of the knot complement (see Conjecture ~\ref{SC}(a)). As far as the authors know, the Slope Conjecture has been proved for knots with up to 10 crossings ~\cite{Gar11b}, adequate knots ~\cite{FKP11}, 2-fusion knots ~\cite{GvdV}, some pretzel knots ~\cite{LV} and a family of Montesinos knots ~\cite{LYL17}. In ~\cite{MT15}, K. Motegi and T. Takata verify the conjecture for graph knots and prove that it is closed under taking connected sums. In ~\cite{KT15}, E. Kalfagianni and A. T. Tran prove the conjecture is closed under taking the $(p,q)$-cable with certain conditions on the colored Jones polynomial, and they formulate the Strong Slope Conjecture (see Conjecture ~\ref{SC}(b)).

In this article we verify the Slope Conjecture and the Strong Slope Conjecture for 3-string Montesinos knots $M([r_0,\cdots,r_m],[s_0,\cdots,s_p],[t_0,\cdots,t_q])$ (see Figure 1) with $m,p,q\geq1$ and certain conditions attached (See conditions C(1) and C(2) in Section 2), where
\[
[r_0,\cdots,r_m]=\frac{1}{r_0-\frac{1}{r_1-\frac{1}{\cdots-\frac{1}{r_m}}}},
\]
and $[s_0,\cdots,s_p]$ and $[t_0,\cdots,t_q]$ are defined similarly. Note that our conventions for Montesinos knots coincide with those of ~\cite{MK06}.

This article is a succeeding work of ~\cite{LV} and ~\cite{LYL17} (however, the results of this article is not strictly the generalization of that of ~\cite{LV} or ~\cite{LYL17} because we can not loosen the restriction $m,p,q\geq 1$ in Lemma 3.4), and the goal of these articles is to provide more data and evidences to the (Strong) Slope Conjecture for Montesinos knots. The reason to choose the family of Montesinos knots is that as a generalization of 2-bridge knots, it is large and representative, and meanwhile well parameterized. Moreover, C. R. S. Lee and R. van der Veen's method~\cite{LV} to deal with the colored Jones polynomial and its degree and Hatcher and Oertel's algorithm \cite{HO89} to determine the incompressible surfaces of  Montesinos knots pave the way for the proof. The strategy of the proof is straight-forward: we first find out the maximal degree of the colored Jones polynomial and then choose the essential surface which matches the degree by the boundary slope and the Euler characteristic provided by the Hatcher-Oertel algorithm.

As we will see, for 3-string Montesinos knots $M([r_0,\cdots,r_m],[s_0,\cdots,s_p],\newline [t_0,\cdots,t_q])$, the increasing of $m$, $p$, $q$ does not cause much complexity, and like the cases in ~\cite{LV} and ~\cite{LYL17}, $r_0$, $s_0$ and $s_0$, particularly the discriminant $\Delta$ (see Theorem ~\ref{Jones slope} and its proof in Section 3) still dominate the maximal degree of the colored Jones polynomial (Theorem ~\ref{Jones slope}) as well as the selection of the essential surface (Theorem ~\ref{boundary slope}). More specifically, when $\Delta<0$, the degree of colored Jones polynomial is matched by a typical type I essential surface; when $\Delta\geq 0$, it is matched by a type II essential surface, but this type II surface generally (when at least one of  $m$,$p$ and $q$ is greater than 1) does not correspond to a Seifert surface while it does in ~\cite{LV} and ~\cite{LYL17}.

\section{The Slope Conjectures}
Let $K$ denote a knot in $S^3$ and $N(K)$ denote its tubular neighbourhood. A surface $S$ properly embedded in the knot exterior $E(K)=S^3-N(K)$ is called \textit{essential} if it is incompressible, $\partial$-incompressible, and non $\partial$-parallel. A fraction $\frac{p}{q}\in \mathbb{Q}\bigcup \{\infty\}$ is a \textit{boundary slope} of $K$ if $pm+ql$ represents the homology class of $\partial S$ in the torus $\partial N(K)$, where $m$ and $l$ are the canonical meridian and longitude basis of $H_1(\partial N(K))$. The \textit{number of sheets} of $S$, denoted by $\sharp S$, is the minimal number of intersections of $\partial S$ and the meridional circle of $\partial N(K)$.

For the colored Jones polynomial, we follow the convention of~\cite{LV}, which denotes the unnormalized $\textit{n-colored Jones polynomial}$ by $J_{K}(n;v)$. See Section 3 for details. Its value on the trivial knot is defined to be $[n]=\frac{v^{2n}-v^{-2n}}{v^2-v^{-2}}$, where $v= A^{-1}$, and A is the variable of the Kauffman bracket. The maximal degree of $J_{K}(n)$ is denoted by $d_{+}J_{K}(n)$.

A significant result made by S. Garoufalidis and T. Q. T. Le shows that the colored Jones polynomial is \textit{q}-holonomic ~\cite{GL05}. Furthermore, the degree of the colored Jones polynomial is a quadratic quasi-polynomial~\cite{Gar11a}, which can be stated as follow.

\begin{thm}~\cite{Gar11a}
For any knot $K$, there exist an integer $p_{K}\in \mathbb{N}$ and quadratic polynomials $Q_{K,1}\ldots Q_{K,p_K} \in \mathbb{Q}[x]$ such that $d_{+} J_{K}(n)=Q_{K,j}(n)$ if $n=j\ (mod\ p_{K})$ for sufficiently large $n$.
\end{thm}

Then the Slope Conjecture and the Strong Slope Conjecture can be formulated as follows:
\begin{conj}\label{SC}
In the context of the above theorem, set $Q_{K,j}=a_{j}x^2+ 2b_{j}x+ c_j$, then for each $j$ there exists an essential surface $S_j \subset S^3-K$, such that:

a.(Slope Conjecture~\cite{Gar11b}) $a_j$ is a boundary slope of $S_j$,

b.(Strong Slope Conjecture~\cite{KT15}) $ b_j=\frac{\chi(S_j)}{\sharp S_j}$, where $\chi(S_j)$ is the Euler characteristic of $S_j$.

\end{conj}

\begin{figure}[!ht]
\centering

\includegraphics{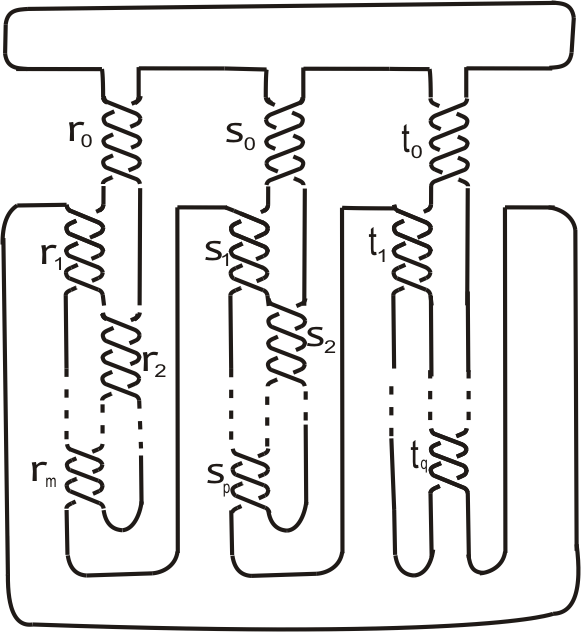}

\caption{The Montesinos Knot $M([r_0,\cdots,r_m],[s_0,\cdots,s_p],[t_0,\cdots,t_q])$}
\end{figure}

A Montesinos knot is a knot formed by putting rational tangles together in a circle (see Figure 1). We denote the Montesinos knot obtained from rational tangles $R_{1}$, $R_{2}$, \ldots $R_{N}$ by $ M(R_{1}, R_{2},..., R_{N})$. For properties about Montesinos knots, the reader can refer to ~\cite{BZ}. It is known that all Montesinos knots are semi-adequate~\cite{LT}, and the Slope Conjecture has been verified for adequate knots~\cite{FKP11}. So we focus on a family of A-adequate and non-B adequate knots $M([r_0,\cdots,r_m],[s_0,\cdots,s_p],[t_0,\cdots,t_q])$ with conditions on $\{r_i\}$, $\{s_j\}$ and $\{t_k\}$ as follows.\\
{\bf C(1)}: For $M$ to be knots rather than links, according to ~\cite{HO89}(Pg.456), let $r_m$, $s_p$ and $t_k$ be odd integers and all the rest be even integers such that $[r_0,\cdots,r_m]$,$[s_0,\cdots,s_p]$ and $[t_0,\cdots,t_q]$ are all of the type $\frac{odd}{odd}$.\\
{\bf C(2)}: For $M$ to be A-adequate, let $s_0$, $t_0$ be positive integers and all the rest be negative integers.

Note that each of the above conditions is sufficient but unnecessary.

Our main theorem is stated as follow.
\begin{thm}
The Slope Conjecture and the Strong Slope Conjecture are true for the Montesinos knots $M([r_0,\cdots,r_m],[s_0,\cdots,s_p],[t_0,\cdots,t_q])$ with  $\{r_i\}$, $\{s_j\}$ and $\{t_k\}$ satisfying conditions C(1) and C(2).
\end{thm}

 This theorem will be proved directly from the following two theorems. The first is about the degree of the colored Jones polynomial and the second is about the essential surface.

\begin{thm}\label{Jones slope}
Let $K=M([r_0,\cdots,r_m],[s_0,\cdots,s_p],[t_0,\cdots,t_q])$ satisfying conditions C(1), C(2), set $A=-\frac{r_0+s_0+2}{2}, B=-(r_0+2), C=-\frac{r_0+t_0+2}{2}, \Delta=4AC-B^2 $.  \\
(1) If $\Delta<0$, then $p_K =\frac{s_0+t_0}{2}$, and
\[
\begin{split}
&d_{+} J_{K}(n)= Q_{K,l}\\
&=[\frac{2t_0^{2}}{s_0+t_0}-2(r_0+t_0+2)+6-2(m+p+q)-2(\sum_{even}^{m}r_i + \sum_{even}^{p}s_j +\sum_{even}^{q}t_k)]n^{2}\\
&+2[r_0 + 2(m+p+q)+ \sum_{i=1}^{m}r_i + \sum_{j=1}^{p}s_j +\sum_{k=1}^{q}t_k]n\\
& - \frac{s_0+t_0}{2}\alpha_{l}^2- (s+t)\alpha_{l} -2(m+p+q)-2-2(\sum_{odd}^{m}r_i + \sum_{odd}^{p}s_j + \sum_{odd}^{q}t_k).
\end{split}
\]
where $\alpha_{l}$ is defined as follows. Let $0\leq l < \frac{s_0+t_0}{2}$ such that  $n= l\ mod\  \frac{s_0+t_0}{2}$ , and set $v_l$ to be the odd number nearest to $\frac{2t_{0}}{s_0+t_0}l$, then we set $\alpha_l = -\frac{2t_0}{s_0+t_0}l+v_{l}-1 $. Note $p_K=\frac{s_0+t_0}{2}$ is a period of $d_{+} J_{K}(n)$ but may not be the least one. And $\sum_{even}^{m}$ ($\sum_{odd}^{m}$) means the summation is over all positive even (odd) numbers not greater than $m$, and $\sum_{even}^{p}$ ($\sum_{odd}^{p}$) and $\sum_{even}^{q}$ ($\sum_{odd}^{q}$) are defined similarly.  \\
(2) If $\Delta\geq0$, then $p_{K}=1$ and
\[
\begin{split}
d_{+}J_{K}(n)=&[6-2(m+p+q)-2(\sum_{even}^{m}r_i + \sum_{even}^{p}s_j+ \sum_{even}^{q}t_k)]n^2\\
&+ 2[2(m+p+q)-4+ (\sum_{i=1}^{m}r_i + \sum_{j=1}^{p}s_j +\sum_{k=1}^{q}t_k)]n \\
&+2-2(m+p+q)-2(\sum_{odd}^{m}r_i + \sum_{odd}^{p}s_j+ \sum_{odd}^{q}t_k).
\end{split}
\]
\end{thm}

\begin{thm}\label{boundary slope}
Under the same assumptions as Theorem 2.4, \\
(1) When $\Delta<0 $, there exists an essential surface $S_1$ with boundary slope\\
 \[bs(S_1)=\frac{2t_0^{2}}{s_0+t_0}-2(r_0+t_0+2)+6-2(m+p+q)-2(\sum_{even}^{m}r_i + \sum_{even}^{p}s_j + \sum_{even}^{q}t_k),\]
 and
 \[\frac{\chi_{(S_1)}}{\sharp S_1}=r_0 + 2(m+p+q)+ \sum_{i=1}^{m}r_i + \sum_{j=1}^{p}s_j +\sum_{k=1}^{q}t_k.\]
(2) When $\Delta\geq0$, there exists an essential surface $S_2$ with boundary slope\\
\[bs(S_2)=6-2(m+p+q)-2(\sum_{even}^{m}r_i + \sum_{even}^{p}s_j+ \sum_{even}^{q}t_k),\]
and
\[\frac{\chi_{(S_2)}}{\sharp S_2}=2(m+p+q)-4+ (\sum_{i=1}^{m}r_i + \sum_{j=1}^{p}s_j +\sum_{k=1}^{q}t_k). \]
\end{thm}

Note that in Theorem ~\ref{Jones slope} the coefficient of the linear term of $d_{+} J_{K}(n)$ is always negative. This actually verifies another conjecture from ~\cite{KT15} for this family of Montesinos knots, which can be stated as follow.
\begin{conj}(Conjecture 5.1, ~\cite{KT15})
In the context of Theorem 2.1 and Conjecture 2.2, for any nontrivial knot in $S^3$, we have $b_{j}\leq 0$.
\end{conj}

\begin{thm}
Conjecture 2.6 is true for the Montesinos knots $M([r_0,\cdots,r_m],\newline [s_0,\cdots,s_p],[t_0,\cdots,t_q])$ satisfying the Condition (1) and (2).
\end{thm}

\section{The Colored Jones Polynomial and Its Maximal Degree}
To compute the colored Jones polynomial of Montesinos knots, Lee and van der Veen introduce the notion of \textit{knotted trivalent graphs} (KTG) in~\cite{LV} (see also~\cite{vdV09,Thu02}). It is a natural generalization of knots and links and makes the skein theory ~\cite{MV} more convenient for Montesinos knots.

\begin{figure}[!ht]
\centering
\includegraphics{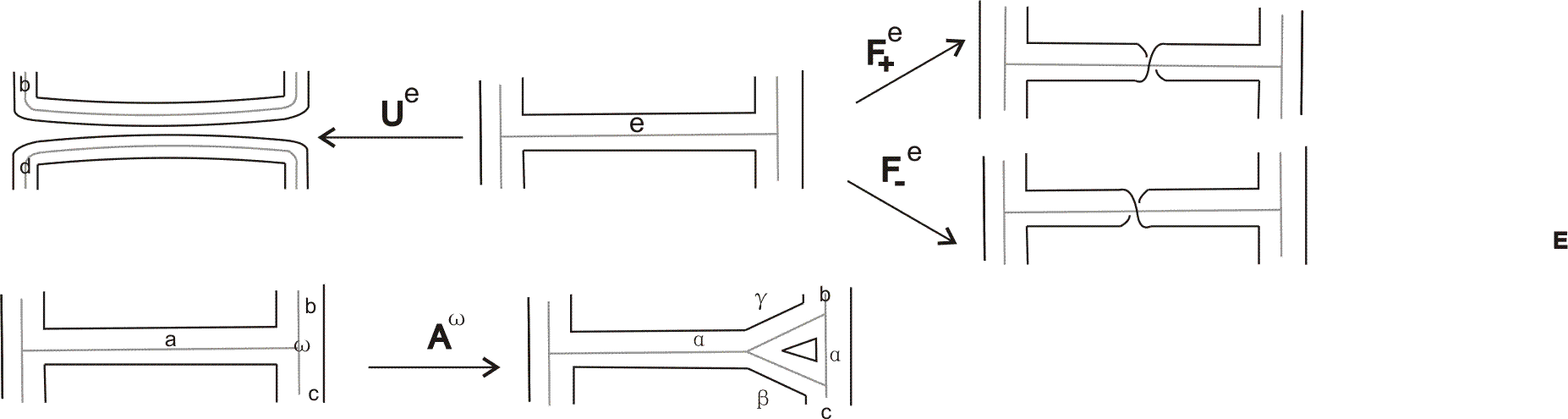}
\caption{Operations on KTG: framing change $F$ and unzip $U$ applied to an edge $e$, triangle move $A^{\omega}$ applied to a vertex $\omega$.}
\end{figure}

\begin{dfn}
~\cite{LV}

(1) A \textit{framed graph} is a one dimensional simplicial complex $\Gamma$ together with an embedding $\Gamma \rightarrow \Sigma$ of $\Gamma$ into a  surface with boundary $\Sigma$ as a spine.

(2) A \textit{coloring} of $\Gamma$ is a map $\sigma: E(\Gamma)\rightarrow\mathbb{N}$, where $E(\Gamma)$ is the set of edges of $\Gamma$.

(3) A \textit{knotted trivalent graph} (KTG) is a trivalent framed graph embedded as a surface into $\mathbb{R}^3$, considered up to isotopy.
\end{dfn}
The advantage of KTGs over knots or links is that they support powerful operations. In this article we will need the following three types of operations, the \textit{framing change} $F^{e}_{\pm}$, the \textit{unzip} $U^e$, and the \textit{triangle move} $A^{\omega}$, as illustrated in Figure 2.

The important thing is that these three types of operations are sufficient to produce any KTG from the $\theta$ graph (see the far left in Figure 3).

\begin{thm}
~\cite{vdV09,Thu02} Any KTG can be generated from the $\theta$ graph by repeatedly applying the operations $F_{\pm}$, $U$ and $A$ defined above.

\end{thm}

By the above theorem, one can define the colored Jones polynomial of any KTG once he or she fixes the value of any colored $\theta$ graph and describes how it varies under the the above operations.
\begin{dfn}

~\cite{LV} The colored Jones polynomial of a KTG $\Gamma $ with coloring $\sigma$, denoted by $\langle\Gamma,\sigma\rangle$, is defined by the equations as follows.
\[
\langle\Theta; a, b, c\rangle= O^{\frac{a+b+c}{2}}\left[
\begin{matrix}
 &\frac{a+b+c}{2}& &\\
\frac{-a+b+c}{2}&\frac{a-b+c}{2}&\frac{a+b-c}{2}&
\end{matrix}
\right],
\]

\[
\langle F_{\pm}^{e}(\Gamma), \sigma\rangle= f(\sigma(e))^{\pm 1}\langle\Gamma, \sigma\rangle,
\]

\[
\langle U^{e}(\Gamma), \sigma\rangle= \langle \Gamma, \sigma\rangle \sum_{\sigma(e)}\frac{O^{\sigma(e)}}{\langle\theta; \sigma(e), \sigma(b), \sigma(d)\rangle},
\]

\[
\langle A^{\omega}(\Gamma), \sigma\rangle= \langle \Gamma, \sigma\rangle \Delta(a, b, c, \alpha, \beta, \gamma).
\]

\end{dfn}
Particularly, a knot $K$ is a $0$-frame KTG without vertices, and the colored Jones polynomial of $K$ is defined to be $J_{K}(n+1)=(-1)^{n}\langle K, n \rangle$, where $n$ is the color of the single edge of $K$, and $(-1)^n$ is to normalize the unknot as $J_{O}(n)=[n]$.

In above formulas, the quantum integer \[[k]=\frac{v^{2k}-v^{-2k}}{v^2-v^{-2}},\ and \ [k]!=[k][k-1]\ldots[1].\] The symmetric multinomial coefficient is defined as:
\[\left[
\begin{matrix}
 a_1+ a_2+ \ldots a_r\\
 a_1,a_2,\ldots, a_r
\end{matrix}
\right]=\frac{[a_{1}+a_{2}+...+a_{r}]!}{[a_{1}]!\ldots[a_{r}]!}.
\]
The value of the $k$-colored unknot is defined as: \[O^k=(-1)^{k}[k+1]=\langle O,k \rangle. \]
The the framing change $f$ is defined as:
\[f(a)=(\sqrt{-1})^{-a} v^{-\frac{1}{2}a(a+2)}.\]
The summation in the equation of unzip is over all admissible colorings of the edge $e$ that has been unzipped.
$\Delta$ is the quotient of the $6j$-symbol and the $\theta$, and \[\Delta(a,b,c,\alpha,\beta,\gamma)=\Sigma\frac{(-1)^z}{(-1)^{\frac{a+b+c}{2}}}\left[
\begin{matrix}
z+1 \\
 \frac{a+b+c}{2}+1
\end{matrix}
\right]\left[\begin{matrix}
 \frac{-a+b+c}{2}\\
z-\frac{a+\beta+\gamma}{2}
\end{matrix}
\right]\left[\begin{matrix}
 \frac{a-b+c}{2}\\
z-\frac{\alpha+b+\gamma}{2}
\end{matrix}
\right]\left[\begin{matrix}
 \frac{a+b-c}{2}\\
z-\frac{\alpha+\beta+c}{2}
\end{matrix}
\right].
\]
The range of the summation in above formula is indicated by the binomials. Note that this $\Delta$ is not the one in Theorem 2.4.

The above definition agrees with the integer normalization in~\cite{Cos}, where F. Costantino shows that $\langle \Gamma, \sigma\rangle$ is a Laurent polynomial in $v$ independent of the choice of operations to produce the KTG.

\begin{figure}[!ht]~\label{Jones2}
\centering

\includegraphics{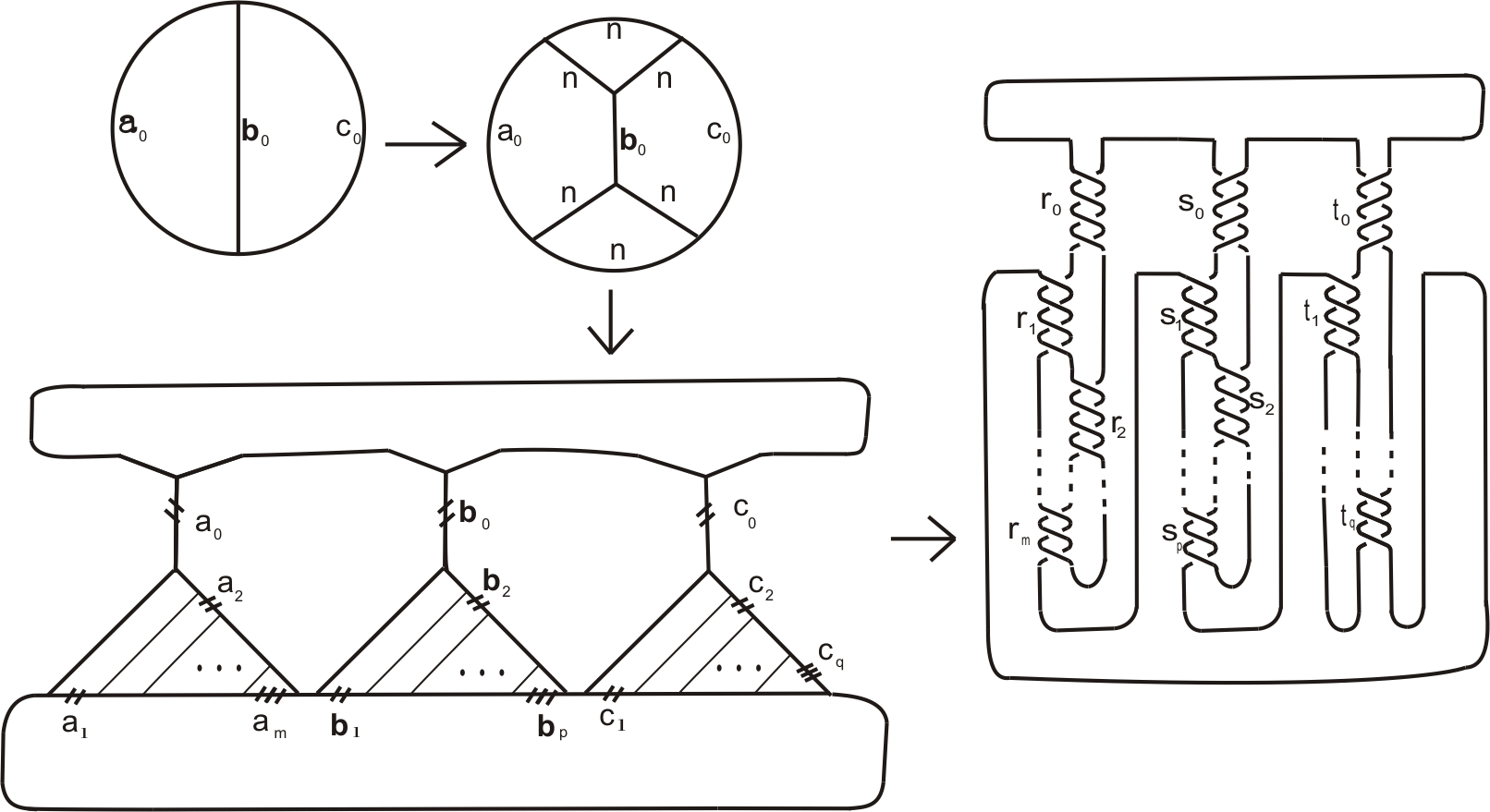}

\caption{The operations to produce the knot $K$ from a $\theta$ graph.  }
\end{figure}

As illustrated in Figure ~\ref{Jones2}, we obtain the colored Jones polynomial of the knot $K=M([r_0,\cdots,r_m],[s_0,\cdots,s_p],[t_0,\cdots,t_q])$ as follows. Starting from a $\theta$ graph , we first apply two $A$ moves, then $(m+p+q)$ $A$ moves on the three vertices of the lower triangle of the second graph, then one $F$ move on each of the edges labelled by $a_i$, $b_j$ and $c_k$, then unzip these twisted edges. The edges without labelling are actually colored by $n$. Note that an unzip applied to a twisted edge produces two twisted bands, each of which has the same twist number of the unzipped edge. Finally, to get the $0$-frame colored Jones polynomial we need to cancel the framing produced by the operations and the writhe of the knot, which are denoted by $F(K)$ and $writhe(K)$ respectively and computed as follows.
\[F(K)=\sum_{i=0}^{m}r_i + \sum_{j=0}^{p}s_j +\sum_{k=0}^{q}t_k\],
\[writhe(K)=\sum_{i=0}^{m}(-1)^{i+1}r_i + \sum_{j=0}^{p}(-1)^{j+1}s_j +\sum_{k=0}^{q}(-1)^{k+1}t_k .\]
so the result should be multiplied by
\[f(n)^{-2F(K)-2writhe(K)}=f(n)^{-4(\sum_{odd}^{m}r_i + \sum_{odd}^{p}s_j +\sum_{odd}^{q}t_k)}.\]

\begin{lem}\label{CJP}

The colored Jones polynomial of the Montesinos knot $K=M([r_0,\cdots,r_m],[s_0,\cdots,s_p],[t_0,\cdots,t_q])$ is

\[
\begin{split}
  &J_{K}(n+1)=\\
  &(-1)^{n}f(n)^{-4(\sum_{odd}r_i + \sum_{odd}s_j+ \sum_{odd}t_k) }\!+\sum_{(a_i, b_j, c_k)\in D_n}^{}
  \langle\theta;\!a_0,\! b_0,\! c_0\rangle\Delta^{2}(a_0,\! b_0,\! c_0,\! n,\! n,\! n)\\
  &\prod_{i=0}^{m-1}\Delta(a_{i},\! n,\! n,\! a_{i+1},\! n,\! n) \prod_{j=0}^{p-1}\Delta(b_{j},\! n,\! n,\! b_{j+1},\! n,\! n) \prod_{k=0}^{q-1}\Delta(c_{k},\! n,\! n,\! c_{k+1},\! n,\! n)\\
  &\prod_{i=0}^{m}f^{r_i}(a_i) \prod_{j=0}^{p}f^{s_j}(b_j) \prod_{k=0}^{q}f^{t_k}(c_k)
   \ \prod_{i=0}^{m}O^{a_i}\langle\theta;\!a_i,\! n,\! n\rangle^{-1}\prod_{j=0}^{p}O^{b_j}\langle\theta;\!b_j,\!n,\!n\rangle^{-1}\\
  &\prod_{k=0}^{q}O^{c_k}\langle\theta;\!c_k,\! n,\! n\rangle^{-1},
\end{split}
\]
where the domain $D_n$ is defined such that $a_i$,$b_j$,$c_k$ are all even with $0\leq a_i, b_j, c_k \leq 2n$, and $a_0$,$b_0$,$c_0$ satisfy the triangle inequality.

\end{lem}

To find out the maximal degree of the colored Jones polynomial, we need to analyze the the factors of the summands. The following lemma and is from~\cite{LV}.

\begin{lem}\label{degree}~\cite{LV}
\[
d_{+}f(a)=-\frac{a(a+2)}{2},
\]
\[
d_{+}O^a = 2a,
\]
\[
d_{+}\langle\theta; a, b, c\rangle= a(1-a)+b(1-b)+c(1-c)+\frac{(a+b+c)^{2}}{2},
\]

\[
\begin{split}
d_{+}\Delta(a, b, c, \alpha, \beta, \gamma)&= g(m+1,\frac{a+b+c}{2}+1)+ g(\frac{-a+b+c}{2},m-\frac{a+\beta+\gamma}{2})\\
&+g(\frac{a-b+c}{2},m-\frac{\alpha+b+\gamma}{2})+g(\frac{a+b-c}{2},m-\frac{\alpha+\beta+c}{2}),\\
\end{split}
\]
where $g(n,k)=2k(n-k)$ and $2m=a+b+c+\alpha+\beta+\gamma-max(a+\alpha,b+\beta,c+\gamma).$

\end{lem}

Now we can apply Lemma 3.4 and 3.5 to prove Theorem 2.4.

\begin{proof}[{\bf Proof of Theorem 2.4}]
Note that the maximal degree of $J_{K}(n+1)$ satisfies the inequality below.
\[d_{+}J_{K}(n+1)\leq max_{(a_i, b_j, c_k,)\in D_{n}}\Phi(a_0,\cdots,a_m,b_0,\cdots,b_p,c_0,\cdots,c_q ),\]
where
\[
\begin{split}
&\Phi(a_0,\cdots,a_m,b_0,\cdots,b_p,c_0,\cdots,c_q )\\
=&-4(\sum_{odd}^{m}r_i + \sum_{odd}^{p}s_j+ \sum_{odd}^{q}t_k)d_+f(n)+d_{+}\langle\theta; a_0, b_0, c_0\rangle+2d_{+}\Delta(a_0, b_0, c_0, n, n, n)\\
    &+\sum_{i=0}^{m-1}d_+\Delta(a_{i}, n, n, a_{i+1}, n, n)+\sum_{j=0}^{p-1}d_+\Delta(b_{j}, n, n, b_{j+1}, n, n)\\
    &+\sum_{k=0}^{q-1}d_+\Delta(c_{k}, n, n, c_{k+1}, n, n)+\sum_{i=0}^{m}r_{i}d_{+}f(a_i)+\sum_{j=0}^{p}s_{j}d_{+}f(b_j)+\sum_{k=0}^{q}t_{k}d_{+}f(c_k)\\
    &+\sum_{i=0}^{m}d_+O^{a_i}+\sum_{j=0}^{p}d_+O^{b_j}+\sum_{k=0}^{q}d_+O^{c_k}-\sum_{i=0}^{m}d_+\langle\theta; a_i, n, n\rangle-\sum_{j=0}^{p}d_+\langle\theta; b_j, n, n\rangle\\
    &-\sum_{k=0}^{q}d_+\langle\theta; c_k, n, n\rangle.
\end{split}
\]
$\Phi(a_0,\cdots,a_m,b_0,\cdots,b_p,c_0,\cdots,c_q )$ is the highest degree of each term of the summation in Lemma ~\ref{CJP}. The equality holds when $\Phi$ has only one maximum or when it has multiple maxima and the coefficients of the maximal degree terms do not cancel out.

Generally, finding $max_{(a_i, b_j, c_k,)\in D_{n}}\Phi(a_0,\cdots,a_m,b_0,\cdots,b_p,c_0,\cdots,c_q )$ is a problem of quadratic integer programming, which is quite a involved topic~\cite{GvdV}. In this case however, it can be solved by observing its monotonicity.

From Lemma ~\ref{degree}, we have

\begin{equation*}
d_+\Delta(a, b, c, n, n, n)=
\begin{cases}
-\frac{1}{2}a^2-a+(a+b+c+2)n-bc & \text{if $a\geq b, c$};\\
-\frac{1}{2}b^2-b+(a+b+c+2)n-ac & \text{if $b\geq a, c$};\\
-\frac{1}{2}c^2-c+(a+b+c+2)n-ab & \text{if $c\geq a, b$},
\end{cases}
\end{equation*}
and
\begin{equation*}
d_+\Delta(a, n, n, b, n, n)=
\begin{cases}
-\frac{a^2}{2}\!-\!a\!-\!b^2\!-\!ab\!+\!2an\!+4bn+\!2n\!-\!2n^2 & \text{if $a\!+\!b\geq 2n$};\\
-\frac{1}{2}b^2+b+2nb & \text{if $a+b\leq 2n$}.\\
\end{cases}.
\end{equation*}

When $m\geq 2$, for $1\leq i\leq m-1$ we have,
\begin{equation*}
\begin{split}
&\partial_{a_i}\Phi(a_0,\cdots,a_m,b_0,\cdots,b_p,c_0,\cdots,c_q )\\
=&a_i\!-\!r_i(a_i\!+\!1\!)\!+\!1\!-\!2n+\!\partial_{a_i}d_+\Delta(a_{i-1},\!n,\!n,a_{i},\!n,\!n)\!+\partial_{a_{i}}d_+\Delta(a_i,\!n,\!n,\!a_{i+1},\!n,\!n)\\
           =&\begin{cases}
             -(r_i+2)a_i-\!r_i-\!a_{i-1}-a_{i+1}+4n & \!\text{if $a_{i-1}\!+\!a_i\geq 2n$ and $a_{i}\!+\!a_{i+1}\geq 2n;$}\\
            -(r_i+1)a_i-a_{i-1}-r_i+2n+\!1 &\! \text{if $a_{i-1}\!+\!a_i\geq 2n$ and $a_{i}\!+\!a_{i+1}\leq 2n;$}\\
            -(r_i+1)a_i-r_i+2n-a_{i+1}+1   &\! \text{if $a_{i-1}\!+\!a_i\leq 2n$ and $a_{i}\!+\!a_{i+1}\geq 2n;$}\\
            -r_i(a_i+1)+2                & \!\text{if $a_{i-1}\!+\!a_i\leq 2n$ and $a_{i}\!+\!a_{i+1}\leq 2n.$}
    \end{cases}
\end{split}
\end{equation*}
For $i=m$ ($m\geq 2$), we have
\begin{equation*}
\begin{split}
 &\partial_{a_m}\Phi(a_0,\cdots,a_m,b_0,\cdots,b_p,c_0,\cdots,c_q )\\
=&a_m-r_m(a_m+1)+1-2n+\partial_{a_m}d_+\Delta(a_{m-1},n,n,a_{m},n,n)\\
=&\begin{cases}
 -(r_m+1)a_m-a_{m-1}-r_m+2n+1  & \text{if $a_{m-1}+a_m\geq 2n$;}\\
 -r_m(a_m+1)+2             & \text{if $a_{m-1}+a_m\leq 2n$.}
 \end{cases}
\end{split}
\end{equation*}
Since $r_i\leq -2$ ($1\leq i\leq m$), from the two equations above it is easy to verify that we have $\partial_{a_i}\Phi>0$ in all cases.

When $m=1$, by a similar calculation we have  $\partial_{a_1}\Phi>0$. So we can conclude that $\partial_{a_i}\Phi(a_0,\cdots,a_m,b_0,\cdots,b_p,c_0,\cdots,c_q )>0$ with $m\geq 1$ and $1\leq i \leq m$.
Similarly, we have $\partial_{b_j}\Phi>0$, $\partial_{c_k}\Phi>0$ with $p,q\geq 1$ and $1\leq j\leq p$, $1\leq k\leq q$. So $\Phi(a_0,\cdots,a_m,b_0,\cdots, b_p,c_0,\cdots,c_q )$ achieves its maxima only when $a_i=2n$, $b_j=2n$ and $c_k=2n$, where $m,p,q\geq 1$, $1\leq i\leq m$, $1\leq j\leq p$ and $1\leq k\leq q$.
Further, we have
\begin{equation*}
\begin{split}
&\partial_{a_0}\Phi(a_0,b_0,c_0,2n,\cdots,2n)\\
=&b_0\!+\!c_0\!+\!2\!-\!2n\!-\!(r_0\!+\!1)a_0\!+\!2\partial_{a_0}d_+\Delta(a_0,\!b_0,\!c_0,\!n,\!n,\!n)+\partial_{a_{0}}d_+\Delta(a_0,\!n,\!n,\!2n,\!n,\!n)\\
=&\begin{cases}
-(r_0+2)(a_0+1)-a_0+b_0+c_0+1 & \text{if $a_0\geq b_0, c_0$};\\
-(r_0+1)(a_0+1)+b_0-c_0+2 & \text{if $b_0\geq a_0, c_0$};\\
-(r_0+1)(a_0+1)+c_0-b_0+2 & \text{if $c_0\geq a_0, b_0$}.
\end{cases}
\end{split}
\end{equation*}

Since $r_0\leq-2$, we always have $\partial_{a_0}\Phi(a_0,b_0,c_0,2n,\cdots,2n)> 0$. See Figure 4. The domain of the real function $\Phi(a_0,b_0,c_0,2n,\cdots,2n)$ is the hexahedron $AB'CD'C'$. Note that for any  $(b'_0,c'_0)$ ($b'_0$ and $c'_0$ are even integers and $0\leq b'_0,c'_0\leq 2n$), there exists an even integer $a'_0$ ($0\leq a'_0\leq 2n$) such that  $(a'_0,b'_0,c'_0)$ is in the triangle region $AB'C$ or $B'CC'$. So $\Phi(a_0,b_0,c_0,2n,\cdots,2n)$ must achieve its maxima in the triangle region $AB'C$ or $B'CC'$.
Note that in the tetrahedron $AB'CC'$ we have
\begin{equation*}
\partial_{b_0}\Phi(a_0,b_0,c_0,2n,\cdots,2n)=a_0-b_0-c_0+1-s_0-s_0 b_0<0.
\end{equation*}

\begin{figure}[!ht]~\label{cube}
\centering
\includegraphics{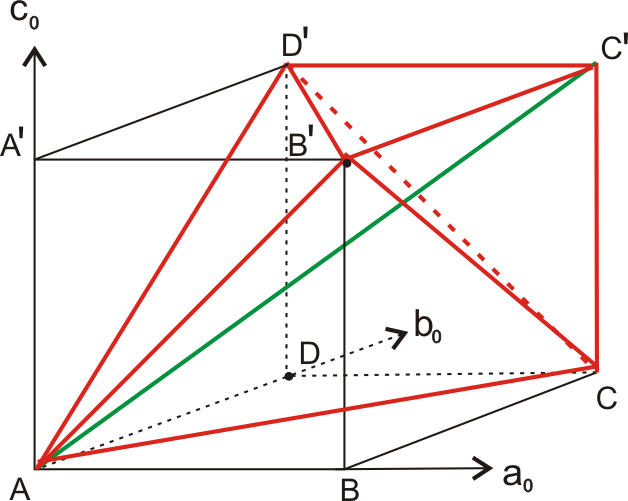}

\caption{The feasible region of real $\Phi(\!a_0,\!b_0,\!c_0,\!2n,\!\cdots,\!2n\!)$ is the hexahedron $AB'CD'C'$ in the cube $ABCD-A'B'C'D'$ with edge length $2n$.}

\end{figure}

So any maximum of $\Phi(a_0,b_0,c_0,2n,\cdots,2n)$ must occur in the triangle region $AB'C$ with $a_0=b_0 + c_0$.
Now we focus on the following 2-variable function $R(b_0,c_0)$ restricted in the triangle domain $T_{n}=\{(b_0,c_0)\mid b_0,c_0\geq 0, b_0+c_0\leq 2n\}$.
\begin{equation}
\begin{split}
&R(b_0, c_0)=\Phi(b_0+c_0, b_0, c_0, 2n,\cdots, 2n)\\
=&-\frac{r_0+s_0+2}{2}b_{0}^{2}-(r_0+2)b_{0}c_0-\frac{r_0+t_0+2}{2}c_0^{2}-(r_0+s_0)b_0-(r_0+t_0)c_0\\
&+[6-2(m+p+q)]n^2 + 4n -2n(n+1)(\sum_{i=1}^{m}r_i + \sum_{j=1}^{p}s_j +\sum_{k=1}^{q}t_k)\\
& +2n(n+2)(\sum_{odd}^{m}r_i + \sum_{odd}^{p}s_j+ \sum_{odd}^{q}t_k).\\
\end{split}
\end{equation}

To analyse $R(b_0, c_0)$, we set $A=-\frac{r_0+s_0+2}{2}$, $B=-(r_0+2)$, $C=-\frac{r_0+t_0+2}{2}$ and $\Delta=4AC-B^2 $.

Although Theorem 2.4 (also Theorem 2.5) is divided into two cases by the range of $\Delta$, it is more natural to divide Theorem 2.4 by the range of $A$ and $C$ in its proof. In fact, if $A=-\frac{r_0+s_0+2}{2}\geq 0$, then
\[
\begin{split}
\Delta=4AC-B^2=(s_0+t_0)(r_0+2)+s_0t_0
      \leq -s_0(s_0+t_0)+s_0t_0=-s_0^2<0
\end{split}
\]
Similarly, when $C\geq 0$ we have $\Delta<0$. So
\[\Delta<0 \Longleftrightarrow A \ \text{or} \ C \geq 0, \text{or}\ A \ \text{and}\ C <0 \ \text{with}\ \Delta <0.\]
In the right side of ``$\Longleftrightarrow$'', the first condition corresponds to the part (1) of this proof, and the second corresponds to (2)a.

{\bf (1)} When $ A\geq0 $ or $ C\geq0 $, we have
\[\partial_{b_0}R=-(r_0+s_0+2)b_0-(r_0+2)c_0-(r_0+s_0)>0\]
 or \[\partial_{c_0}R=-(r_0+t_0+2)b_0-(r_0+2)c_0-(r_0+t_0)>0\] respectively.

Then the maxima must be on the line $b_0+c_0=2n$. Set $Q(b)=R(b_0,2n-b_0)$, we have

\[
\begin{split}
\ Q(b_0)&=R(b_0,2n-b_0)\\
&=-\frac{s_0+t_0}{2}b_0^{2}+[2t_0n-s_0+t_0]b_0-2(r_0+t_0+2)n^2-2(r_0+t_0)n\\
&+[6-2(m+p+q)]n^2 + 4n -2n(n+1)(\sum_{i=1}^{m}r_i + \sum_{j=1}^{p}s_j +\sum_{k=1}^{q}t_k)\\
& +2n(n+2)(\sum_{odd}^{m}r_i + \sum_{odd}^{p}s_j+ \sum_{odd}^{q}t_k).
\end{split}
\]

$Q(b_0)$ is a quadratic function in $b_0$ with negative leading coefficient, and its real maximum is at $\hat{b_{0}}=\frac{2t_0n-s_0+t_0}{s_0+t_0}$, $\hat{b_{0}}\in (0,2n)$ for $n$ sufficiently large.
Since we need the formula for $d_{+}J_{K}(n)$ rather than $d_{+}J_{K}(n+1)$, we set $n+1=N$, and let $N=h(\frac{s_0+t_0}{2})+l$, where $0\leq l< \frac{s_0+t_0}{2}$, then $\hat{b_{0}}= t_0 h-1+\frac{2t_0 l}{s_0+t_0}$.
Let $\bar{b_{0}}$ be the even number nearest to $\hat{b_{0}}$, then we have $\bar{b_{0}}= t_0 h-1+v_{l}$, where $v_{l}$ is the odd number nearest to $\frac{2t_{0}}{s_0+t_0}l$ ($v_l$ is just the source of the periodicity in the case (1) of Theorem 2.4),
 so $\bar{b_{0}}=\frac{2t_0}{s_0+t_0}N-\frac{2t_0}{s_0+t_0}l+v_{l}-1$. Set $\bar{b_{0}}=\frac{2t_0}{s_0+t_0}N + \alpha_l $, where $\alpha_l = -\frac{2t_0}{s_0+t_0}l+v_{l}-1 $.Then we have
\[
\begin{split}
& max_{(a_i,b_j,c_k)\in D_{n}}\Phi(a_i,b_j,c_k)= Q(\bar{b_{0}})\\
=&[\frac{2t_0^{2}}{s_0+t_0}-2(r_0+t_0+2)+6-2(m+p+q)-2(\sum_{even}^{m}r_i + \sum_{even}^{p}s_j \sum_{even}^{q}t_k)]N^{2}\\
&+2[r_0 + 2(m+p+q)+ \sum_{i=1}^{m}r_i + \sum_{j=1}^{p}s_j +\sum_{k=1}^{q}t_k]N\\
& - \frac{s_0+t_0}{2}\alpha_{l}^2- (s+t)\alpha_{l} -2(m+p+q)-2-2(\sum_{odd}^{m}r_i + \sum_{odd}^{p}s_j \sum_{odd}^{q}t_k).
\end{split}
\]

When $\hat{b_{0}}$ is not odd, the maximum is unique. Otherwise, $\Phi$ has exactly 2 maxima, we need to consider the possibility that the coefficients of the 2 maximal-degree terms may cancel out. From Lemma 3.4 and Definition 3.3, and the fact that we take $a_i =b_j=c_k=2n$ when $ i, j, k\geq 1$, it is easy to see that for the leading coefficient of each term of the summation, without counting the factors independent of $a_0$, $b_0$, $c_0$, the $f$'s contribute $(-1)^{\frac{1}{2} (a_0r_0+b_0s_0+c_0t_0)}$, the $\Delta$'s contribute $(-1)^{-\frac{1}{2}(a_0+b_0+c_0)}$, the $O$'s and the $\theta$'s contribute none, and altogether it is
\[C=(-1)^{\frac{1}{2}[(r_0-1)a_0 + (s_0-1)b_0 + (t_0-1)c_0]}.\]
Furthermore, since any maximum of $Q(b_0)$ must occur on $a_0=b_0+c_0=2n$, we have
\[
\begin{split}
\tilde{C}(b_0)&=\frac{1}{2}[(r_0-1)a_0 + (s_0-1)b_0 + (t_0-1)c_0]\\
                &= \frac{1}{2}[2n(r_0 -1)+(s_0-1)b_0 + (t_0-1)(2n-b_0)].
\end{split}
\]
If there are two maxima $Q(b_0)$ and $Q(b_0 +2)$, we must have
\[\tilde{C}(b_0)- \tilde{C}(b_0 +2)= t_0 -s_0. \]
Since $s_0$ and $t_0$ are even, $t_0 -s_0$ must be even and the coefficients of the two maximal terms will not cancel out. So we have $d_{+}J_{K}(n)=Q(\bar{b_{0}})$.

{\bf (2)} When $A<0$ and $C<0$, for any fixed $\tilde{c_0}$, by equation 3.1, $R(b_0,\tilde{c_0})$ is a quadratic function in $b_0$ with negative leading coefficient, whose axis of symmetry (in the plane $c_0=\tilde{c_0}$ of the $b_0 c_0 R$-coordinates) intersects the line $\partial_{b_0}R=0$ and is perpendicular to the $b_0c_0$-plane.

If $r_0<-2$ (the case when $r_0=-2$ will be analysed at the end of this proof), we consider the real value of $R(b_0,c_0)$ on the line $\partial_{b_0}R=0$:
\begin{equation}
  \left\{
   \begin{array}{l}
  R(b_0,c_0)=-\frac{r_0\!+\!s_0\!+\!2}{2}b_0^{2}-\!(r_0\!+\!2)b_0c_0-\!\frac{r_0+\!t_0+\!2}{2}c_0^2\!-\!(r_0\!+\!s_0)b_0\!-\!(r_0\!+\!t_0)c_0+$const$  \\
  \partial_{b_0}R=-(r_0+s_0+2)b_0-(r_0+2)c_0-(r_0+s_0) =0.\\
   \end{array}
   \right.
  \end{equation}
Then we have
\begin{equation}
 \begin{split}
&R\mid _{\partial _{b_0} R=0}= R(b_0,-\frac{r_0+s_0+2}{r_0+2}b_0-\frac{r_0+s_0}{r_0+2})\\
 =&[-\frac{r_0\!+s_0\!+\!2}{2(r_0\!+\!2)^{2}}\Delta]b_0^{2}\!+\!\frac{r_0\!+\!s_0\!+\!2}{(r_0\!+\!2)^{2}}[(r_0\!+\!2)(r_0\!+\!t_0)\!-\!(r_0\!+\!t_0\!+\!2)(r_0\!+\!s_0)]b_0+ const.
 \end{split}
\end{equation}
If we imagine $R(b_0,c_0)$ as the surface of a mountain, then $R\mid _{\partial _{b_0} R=0}$ is just the ridge of it.

If $\Delta \neq 0$, $R\mid _{\partial _{b_0} R=0}$ is a quadratic function in $b$ whose axis of symmetry is perpendicular to the $b_0c_0$-plane at the point $P$ with coordinates
\[(b_0\!\mid_P,c_0\!\mid_P)=(\frac{(\!r_0\!+\!2)(\!r_0\!+\!t_0\!)\!-\!(\!r_0\!+\!t_0+2\!)(\!r_0\!+\!s_0\!)}{\Delta},\frac{(\!r_0\!+\!2\!)(\!r_0\!+\!s_0\!)\!-\!(\!r_0\!+\!s_0\!+\!2\!)(\!r_0\!+\!t_0\!)}{\Delta}).\]($P$ is actually the intersection of $\partial_{b_0}R=0$ and $\partial_{c_0}R=0$).

\begin{figure}[!ht]
\centering

\includegraphics{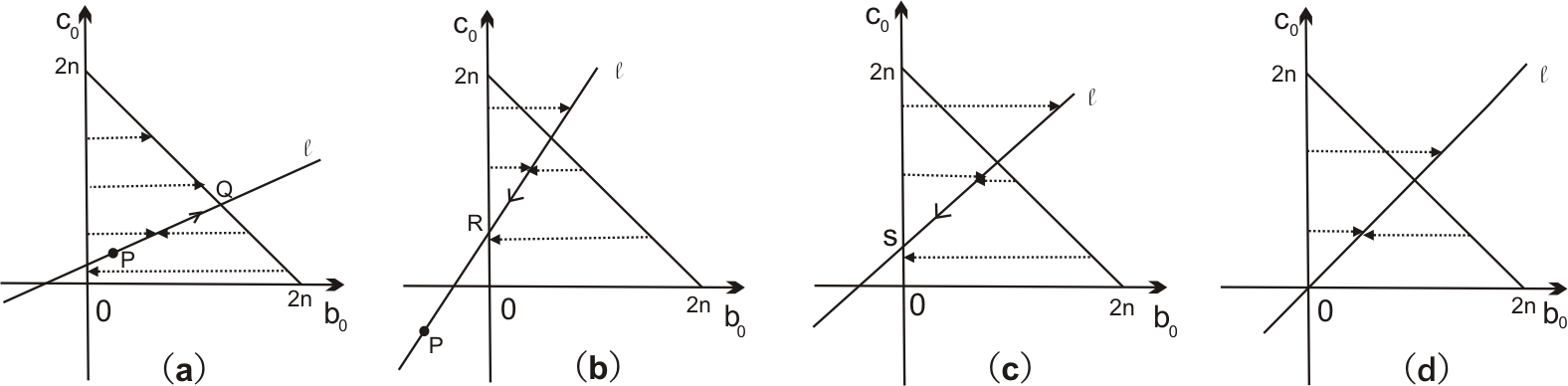}

\caption{$R(b_0,c_0)$ is restricted to the triangle domain $T_n$: $(0,0)$-$(0,2n)$-$(2n,0)$, and the arrows indicate its increasing direction; $\ell$ denotes the line $\partial_{b_0}R=0$.  }
\end{figure}

(a) If $A$ and $C<0, \Delta<0$, by Equation 3.3, $R\!\mid _{\partial _{b_0} R=0}$ is a quadratic function in $b_0$ with positive leading coefficient. And we have $b_0\!\mid_P\geq0$, $c_0\!\mid_P\geq 0$. See Figure 5(a). The arrows indicate the increasing direction of $R(b_0,c_0)$. For sufficiently large $n$, any maximum must be on the segment $[Q,(0,2n)]$ in the line $b_0+c_0=2n$, then the argument will be the same as that of case (1).

(b) If $A,C<0$, and $\Delta>0$, $R\!\mid _{\partial _{b_0} R=0}$ is a quadratic function in $b_0$ with negative leading coefficient. And we have $b_0\!\mid_P\leq 0$, $c_0\!\mid_P\leq 0$. See Figure 5(b). Any maximum must occur on $OR $. Since
 \[R(0,c_0)=-\frac{r_0+t_0+2}{2}c_0^{2}-(r_0+t_0)c_0+ const\]
 and $C=-\frac{r_0+t_0+2}{2}<0$, it is easy to verify that $R(0,c_0)$ decreases in $[0,+\infty)$. And the maximum is unique and must occur at $O=(0,0)$, so
\[
\begin{split}
&d_{+}J_{K}(n+1)=R(0,0)\\
&=[6-2(m+p+q)]n^2 + 4n -2n(n+1)(\sum_{i=1}^{m}r_i + \sum_{j=1}^{p}s_j +\sum_{k=1}^{q}t_k)\\
&+2n(n+2)(\sum_{odd}^{m}r_i + \sum_{odd}^{p}s_j+ \sum_{odd}^{q}t_k).
\end{split}
\]
Let $N=n+1$, we have
\[
\begin{split}
d_{+}J_{K}(N)=&[6-2(m+p+q)-2(\sum_{even}^{m}r_i + \sum_{even}^{p}s_j+ \sum_{even}^{q}t_k)]N^2\\
&+ 2[2(m+p+q)-4+ (\sum_{i=1}^{m}r_i + \sum_{j=1}^{p}s_j +\sum_{k=1}^{q}t_k)]N \\
&+2-2(m+p+q)-2(\sum_{odd}^{m}r_i + \sum_{odd}^{p}s_j+ \sum_{odd}^{q}t_k).
\end{split}
\]

(c) If $A,C<0$ ,$\Delta = 0$, and $(r_0+s_0)^2+(r_0+t_0)^2\neq 0$, $R\mid _{\partial _{b_0} R=0}$ is a decreasing linear function in $b_0$. See Figure 5(c). Any maximum must occur on OS. Since $R(0,c_0)$ decreases in $[0,+\infty)$, the maximum is unique and must be on $O=(0,0)$. So in this case we still have
\[ d_{+}J_{K}(n)=R(0,0).\]

(d) If $A,C<0$, $\Delta =0$, and $(r_0+s_0)^{2}+(r_0+t_0)^{2}=0$, then we immediately have $r_0=-4, s_0=4, t_0=4$, $R(b_0,c_0)=-(b_0-c_0)^{2}+ const$, the maxima are $R(0,0)=R(2,2)=\ldots =R(k,k)$, where $k=n$ when $n$ is even, $k=n-1$ when $n$ is odd. See Figure 5(d). By a similar argument with the end of (1), we can conclude that there are no cancellations between the the highest-degree coefficients, so
\[ d_{+}J_{K}(n)= R(0,0).\]

If $r_0=-2$ (then we must have $A<0$ and $C<0$ and $\Delta>0$, see also Remark 4.4), Equation 3.2 is converted to 
\begin{equation*}
  \left\{
   \begin{array}{l}
  R(b_0,c_0)=-\frac{1}{2}s_0b_0^{2}-\frac{1}{2}t_0c_0^2-(s_0-2)b_0-(t_0-2)c_0+const  \\
  \partial_{b_0}R=-s_0b_0-(s_0-2) =0.\\
   \end{array}
   \right.
  \end{equation*}
Since $\partial_{b_0}R(b_0,c_0)=0 \ \Rightarrow\  b_0=\frac{2-s_0}{s_0}$ and $b_0\leq 0$ when $s_0\geq2$, $R(b_0,c_0)$ must achieve any of its maximum in $c_0$-axis. And $R(0,c_0)=-\frac{1}{2}t_0c_0^2-(t_0-2)c_0+const$ decreases in $[0,+\infty)$. So the unique maximum must be on $O=(0,0)$, and
\[ d_{+}J_{K}(n)= R(0,0).\]

\end{proof}

\section{Boundary Slope and Euler Characteristic  }

 The Hatcher-Oertel edgepath system is actually based on the work of~\cite{HT85}. Like many other topics in geometric topology, the main ideal of the algorithm is to treat the object of study combinatorially. In this mechanism, properly embedded surfaces in a Montesinos knot complement are formed by saddles, and the edgepath system describes how these saddles are combined. For details please refer to ~\cite{HO89, IM07}. Briefly speaking, the \textit{candidate surfaces}, which are the surfaces listed out to include all the essential surfaces with non-empty and non-meridional boundary, are associated to \textit{admissible edgepath systems} in a 1-dimensional diagram $\mathcal{D}$ in the $uv$-plane (see Figure 6(b)). The vertices of $\mathcal{D}$ correspond to projective curve systems $[a, b, c]$ on the 4-punctured sphere carried by the train track in Figure 6(a) via $u=\frac{b}{a+b}$, $v=\frac{c}{a+b}$.
\begin{figure}[!ht]
\centering

\includegraphics{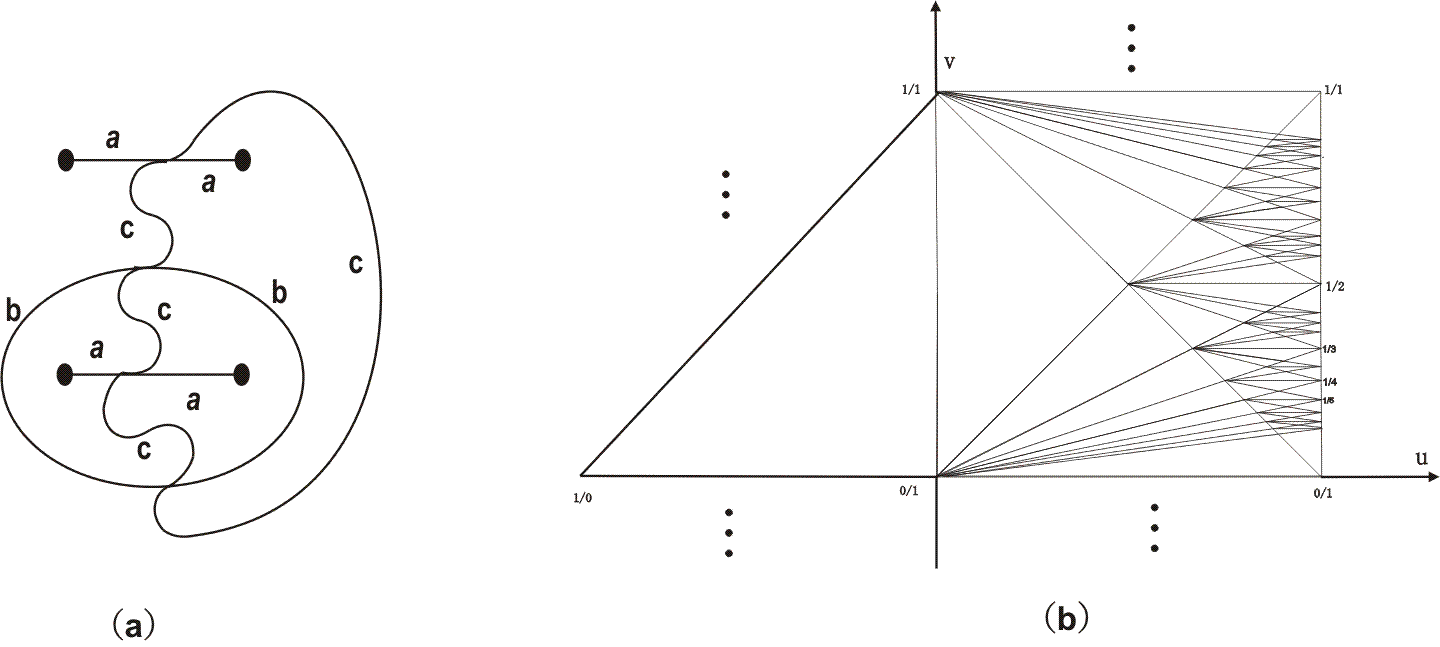}

\caption{(a)The train track in a 4-punctured sphere. (b)The diagram $\mathcal{D}$ in the $uv$-plane. }
\end{figure}

Specifically, the vertices of $\mathcal{D}$ are:

   (1) the vertices corresponding to the arcs with slope $\frac{p}{q}$ denoted by $\langle\frac{p}{q}\rangle$, with the projective curve systems $[1, q-1, p]$, and the \textit{uv}-coordinates $(\frac{q-1}{q}, \frac{p}{q})$,

   (2) the vertices corresponding to the circles with slope $\frac{p}{q}$ denoted by $\langle\frac{p}{q}\rangle ^{\circ}$, with the projective curve systems $[0, p, q]$ and the \textit{uv}-coordinates $(1, \frac{p}{q})$,

   (3) the vertices corresponding to the arcs with slope $\infty$ denoted by $\langle \infty \rangle$, with the \textit{uv}-coordinates $(-1,0)$.

   The edges of $\mathcal{D}$ are:

   (1) the \textit{non-horizontal } edges connecting the vertex $\frac{p}{q}$ to the vertex $\frac{r}{s}$ with $\mid ps- qr\mid =1$, denoted by $\langle \frac{r}{s}\rangle\text{---}\langle\frac{p}{q}\rangle$,

   (2) the \textit{horizontal } edges connecting $\langle\frac{p}{q}\rangle ^{\circ}$ to $\langle\frac{p}{q}\rangle $, denoted by $\langle \frac{p}{q}\rangle \text{---}\langle\frac{p}{q}\rangle ^{\circ} $,

   (3) the \textit{vertical } edges connecting $\langle z\rangle$ to $\langle z\pm 1\rangle$ , denoted by $\langle z\pm 1\rangle\text{---}\langle z\rangle$, here $z\in $ $\mathbb{Z}$,

   (4) the \textit{infinity} edges connecting $\langle z\rangle$ to $\langle \infty \rangle$ denoted by $\langle\infty\rangle\text{---}\langle z\rangle$,

   (5) the \textit{constant} edges which are  points on the horizontal edge $\langle \frac{p}{q}\rangle \text{---}\langle\frac{p}{q}\rangle ^{\circ} $ with the form $ \frac{k}{m} \langle \frac{p}{q}\rangle + \frac{m-k}{m}\langle\frac{p}{q}\rangle ^{\circ} $,

   (6) the \textit{partial} edges which are parts of non-horizontal edges $\langle \frac{r}{s}\rangle\text{---}\langle \frac{p}{q}\rangle $ with the form $\frac{k}{m}\langle \frac{r}{s}\rangle + \frac{m-k}{m} \langle \frac{p}{q}\rangle\text{---} \langle \frac{p}{q}\rangle $.

   An \textit{edgepath} denoted by $\Gamma$ in $\mathcal{D}$ is a piecewise linear path beginning and ending at rational points of $\mathcal{D}$. A \textit {admissible edgepath system} denoted by $\Gamma =(\gamma_{1},\gamma_{2},\ldots, \gamma_{n})$ is an \textit{n}-tuple of edgepaths with the following properties.\\
(E1) The starting point of $\gamma_i$ is on the horizontal edge $\langle \frac{p_i}{q_i} \rangle\text{---}\langle \frac{p_i}{q_i} \rangle^{\circ}$, and if it is not the vertex $\langle \frac{p_i}{q_i}\rangle$, $\gamma_i$ is constant.\\
(E2) $\gamma_i$ is minimal, that is, it never stops or retraces itself, nor does it ever go along two sides of the same triangle of $\mathcal{D}$ in succession.\\
(E3) The ending points of $\gamma_i$'s are rational points of $\mathcal{D}$ with their $\textit{u}$-coordinates  equal and $\textit{v}$-coordinates adding up to zero.\\
(E4) $\gamma_i$ proceeds monotonically from right to left, ``monotonically'' in the weak sense that motion along vertical edges is permitted.

A admissible edgepath system and the corresponding candidate surfaces are called \textit{type I}, \textit{type II} or \textit{type III}, if the $u$-coordinate of the ending points of the admissible edgepath system is positive, zero or negative respectively.

In~\cite{HO89}, a series of candidate surfaces are associated to each admissible edgepath system, then every essential surface in knot complement with non-empty boundary of finite slope is isotopic to one of the candidate surfaces. However, not all of these candidate surfaces are essential, to rule out the inessential surfaces, Hatcher and Oertel developed the notion of \textit{r}-value in ~\cite{HO89}.

\begin{dfn}
The \textit{r}-value of an edge $\langle\frac{t}{s}\rangle\text{---}\langle\frac{p}{q}\rangle$ is defined to be $q-s$, where $0<s<q$. Particularly, when the edge is vertical its \textit{r}-value is $0$. The \textit{r}-value of an edgepath is defined to be the \textit{r}-value of its final edge.
\end{dfn}

By Corollary 2.4 through Proposition 2.10 of ~\cite{HO89}, a series of criteria for incompressibility are established. Here we just extract the part useful for us from Corollary 2.4, Proposition 2.7 and Proposition 2.8(a) of ~\cite{HO89}.
\begin{thm}\label{criteria}~\cite{HO89}
For a 3-string Montesinos knot, a candidate surface associated to the edgepath system \{$\gamma_i$\} is incompressible if it satisfies one of the conditions below:\\
(1) If \{$\gamma_i$\} contains no vertical edges, the cycle of \textit{r}-values of \{$\gamma_i$\} is not the type $(1,1,r_3)$ or $(1,2,r_3)$;\\
(2) If the cycle of \textit{r}-values is the type of $(1,2,r_3)$, the final edges of each edgepath must be all increasing or decreasing;\\
(3) If \{$\gamma_i$\} contains some vertical edges, the cycle of \textit{r}-values is not the type $(0,2,1)$ or $(0,2,0)$.
\end{thm}

The boundary slope of an essential surface $S$ is computed by $\tau(S)- \tau(S_0)$, where $\tau(S)$ is the \textit{total number of twist} (or twist for short) of $S$, and $S_0$ is a Seifert surface in the list of candidate surfaces. For a candidate surface $S$ associated to a admissible edgepath system $\Gamma$, we have~\cite{IM07}
\begin{equation}\label{twist}
  \tau(S)=\sum_{\gamma_{i}\in \Gamma_{non-const}} \sum_{e_{i,j}\in \gamma_{i}}-2\sigma(e_{i,j})|e_{i,j}|.
\end{equation}
 In above formula, $|e|$ is the \textit{length} of an edge $e$, which is defined to be $0$, $1$, or $\frac{k}{m}$ for a constant edge, a complete edge or a partial edge $\frac{k}{m}\langle \frac{r}{s}\rangle + \frac{m-k}{m} \langle \frac{p}{q}\rangle\text{---}\langle \frac{p}{q}\rangle $, respectively. And $\sigma(e)$ is the \textit{sign} of a non-constant edge $e$, which is defined to be $+1$ or $-1$ according to whether the edge is increasing or decreasing (from right to left in \textit{uv}-plane) respectively for a non-$\infty$ edge; for an $\infty$ edge the sign is defined to be $0$.

 For the Euler characteristic of a candidate surface $S$, we use the formulas (3.4) and (3.5) of ~\cite{IM07} to compute the ratio $\frac{-\chi(S)}{\sharp S}$.
 \begin{lem} ~\cite{IM07}\label{Euler}
 Let $S$ be a candidate surface associated to a admissible edgepath system $\Gamma=(\gamma_1,\cdots,\gamma_N)$.\newline
 (1) If $S$ is type I, denote the $u$-coordinate of its ending points by $u_0$, then
 \[\frac{-\chi(S)}{\sharp S}=\sum_{\gamma_{i}\in \Gamma_{non-const}}|\gamma_{i}|+N_{const}-N+(N-2-\sum_{\gamma_{i}\in \Gamma_{const}}\frac{1}{q_{i}})\frac{1}{1-u_{0}}.\]
 Here $\Gamma_{non-const}$ ($\Gamma_{const}$) denotes the set consist of  non-constant (constant) edges, $|\gamma_{i}|$ denotes the length of the edgepath $\gamma_{i}$, $N_{const}$ denotes the number of constant edges in $\Gamma$.\\
 (2) If $S$ is type II, then
 \[\frac{-\chi(S)}{\sharp S}=\sum_{i=1}^{N}(\mid\gamma_{i,>0}\mid) + \mid\Gamma(+0)\mid -2.\]
 Here $\gamma_{i,>0}$ denotes the part of $\gamma$ with positive $u$-coordinate, and $\Gamma(+0)$ denotes the sum of the $v$-coordinates of $\gamma_i$'s when they first reach the $v$-axis.
 \end{lem}

Now we are ready to prove Theorem 2.5.

\begin{proof}[{\bf Proof Theorem 2.5}]
First we note that for $[r_0,\cdots,r_m]$ satisfying the condition C(2), in the diagram $\mathcal{D}$, $\langle[r_0,\ldots, r_{m}]\rangle$ is connected to $\langle[r_0,\ldots, r_{m-1}]\rangle$ by an increasing edge and is connected to $\langle[r_0,\ldots, r_{m}+1]\rangle$ by a decreasing edge from right to left. This is easy to proof by induction and similar facts exist for $[s_0,\cdots,s_p]$ and $[t_0,\cdots,t_q]$.

By the method of~\cite{HO89} (Pg.461),  since the condition C(1) in Section 2 implies that the three tangles of the knot $M([r_0,\cdots,r_m],[s_0,\cdots,s_p],[t_0,\cdots,t_q])$ are all of the form $\frac{odd}{odd}$, we directly find the edgepath system of a Seifert surface $S_0$ as follows. For simplicity, we just use $[r_0,\ldots, r_{m}]$ to denote the vertices instead of $\langle[r_0,\ldots, r_{m}]\rangle$. The edgepaths go from right to left in the same row and the far left vertex of a row is connected to the far right vetex of the row next below. The arrow $\leftharpoonup$ / $\leftharpoondown$ indicates the edge is increasing/ decreasing from right to left.

\[
\begin{split}
\delta_{1}:&\leftharpoonup[r_0, r_1,\ldots, r_{m-1}, -1]\leftharpoondown[r_0, r_1,\ldots, r_{m-1}, -2]\leftharpoondown \cdots \leftharpoondown[r_0, r_1,\ldots,r_{m-1}, r_{m}] \\
           &\leftharpoonup[r_0, r_1,\ldots, r_{m-3}, -1] \leftharpoondown[r_0, r_1,\ldots, r_{m-3}, -2] \leftharpoondown \cdots \leftharpoondown[r_0, r_1,\ldots,r_{m-3}, r_{m-2}]\\
       &\cdots\cdots \\
       &\leftharpoonup[r_0,r_1,r_2,-1]\leftharpoondown[r_0,r_1,r_2,-2]\cdots  \leftharpoondown[r_0,r_1,r_2,r_3]\\
       &\langle 0\rangle\leftharpoonup[r_0,-1]\leftharpoondown[r_0,-2],\ldots, \leftharpoondown[r_0,r_1].\\
\end{split}
\]
\[
\begin{split}
\delta_{2}:&\leftharpoonup[s_0, s_1,\ldots, s_{p-1}, -1]\leftharpoondown[s_0, s_1,\ldots, s_{p-1}, -2]\leftharpoondown \cdots \leftharpoondown[s_0, s_1,\ldots,s_{p-1}, s_{p}] \\
           &\leftharpoonup[s_0, s_1,\ldots, s_{p-3}, -1] \leftharpoondown[s_0, s_1,\ldots, s_{p-3}, -2] \leftharpoondown \cdots \leftharpoondown[s_0, s_1,\ldots,s_{p-3}, s_{p-2}]\\
       &\cdots\cdots \\
       &\leftharpoonup[s_0,s_1,s_2,-1]\leftharpoondown[s_0,s_1,s_2,-2]\cdots  \leftharpoondown[s_0,s_1,s_2,s_3]\\
       &\langle 0\rangle\leftharpoondown[s_0,-1]\leftharpoondown[s_0,-2],\ldots, \leftharpoondown[s_0,s_1].\\
\end{split}
\]
\[
\begin{split}
\delta_{3}:&\leftharpoonup[t_0, t_1,\ldots,t_{q-1}, t_{q}] \\
           &\leftharpoonup[t_0, t_1,\ldots, t_{q-3}, -1] \leftharpoondown[t_0, t_1,\ldots, t_{q-3}, -2] \leftharpoondown \cdots \leftharpoondown[t_0, t_1,\ldots,t_{q-3}, t_{q-2}]\\
       &\cdots\cdots \\
       &\leftharpoonup[t_0,t_1,t_2,-1]\leftharpoondown[t_0,t_1,t_2,-2]\cdots  \leftharpoondown[t_0,t_1,t_2,t_3]\\
       &\langle 0\rangle\leftharpoondown[t_0,-1]\leftharpoondown[t_0,-2],\ldots, \leftharpoondown[t_0,t_1].\\
       \end{split}
\]
Note that the edgepaths of the Seifert surface should avoid the vertices with even denominators, in this case we let $m$ and $p$ be odd and $q$ be even, but the parity of them won't affect our expressions of further results.

The cycle of \textit{r}-value of above edgepath system is $(-r_0-2, s_0, t_0)$. Note that $s_0, t_0\geq 2$, $r_0\leq -2$. Any candidate surface associated to the above edgepath system is essential by Theorem ~\ref{criteria}(1) when $r_0\neq -2$ or by Theorem ~\ref{criteria}(3) when $r_0=-2$.

\begin{rem}
If $r_0=-2$, then $[r_0,-1]=-1$ and there is a vertical edge $\langle 0\rangle\!\leftharpoonup\!\langle -1\rangle$ in the above edgepath $\delta_{1}$ and in the edgepath $\beta_1$ presented in the second part of this proof. Note that when $r_0=-2$, we must have $A=-\frac{1}{2}s_0<0$, $C=-\frac{1}{2}t_0<0$ and $\Delta=\frac{1}{4}s_0t_0>0$. So this situation can only happen in Case (2) of the Theorem 2.5 and won't affect the expressions of our results.
\end{rem}

By formula ~\ref{twist}, the twist of $S_0$ is
\[
\begin{split}
\tau(S_0)&=2[(-r_m-1)+(-1)+(-r_{m-2}-1)+(-1)+\cdots+(-r_1-1)-1]\\
          &+2[(-s_p-1)+(-1)+(-s_{p-2}-1)+(-1)+\cdots+(-s_1-1)+1]\\
          &+2[(-1)+(-t_{q-1}-1)+(-1)+\cdots+(-t_1-1)+1]\\
         &=2-2(m+p+q)-2(\sum_{odd}^{m}r_i + \sum_{odd}^{p}s_j+ \sum_{odd}^{q}t_k).
\end{split}
\]

{\bf (1)}  When $\Delta<0$, we claim that there exists an admissible edgepath system having ending points with $\textit{u}$-coordinate  $u_{0}=\frac{s_0t_0}{s_0t_0+s_0+t_0}$ in $\textit{uv}$-plane. In fact, $u_0$ is just the solution of the equation $v_{1}(u)+v_{2}(u)+ v_{3}(u)=0$, where the linear functions $v=v_{1}(u)$, $v=v_{2}(u)$ and $v=v_{3}(u)$ are determined by the lines through the edges $[\langle -1\rangle\text{--}\langle -\frac{1}{2}\rangle\text{--},...,\text{--}\langle\frac{1}{r_0+1}\rangle]$, $\langle0\rangle\text{---}\langle\frac{1}{s_0+1}\rangle$ and$\langle0\rangle\text{---}\langle\frac{1}{t_0+1}\rangle$, respectively. Denote by $u_{[r_0+1]}$, $u_{[s_0+1]}$ and $u_{[t_0+1]}$ the $\textit{u}$-coordinates of $\langle \frac{1}{r_0+1}\rangle$, $\langle \frac{1}{s_0+1}\rangle$ and $\langle \frac{1}{t_0+1}\rangle$ respectively. With direct calculations we have
\[u_0-u_{[r_0+1]}=\frac{-\Delta}{(r_0+1)(s_0t_0+s_0+t_0)}<0,\]
\[u_0-u_{[s_0+1]}=\frac{-s_0^2}{(s_0+1)(s_0t_0+s_0+t_0)}<0,\]
\[u_0-u_{[t_0+1]}=\frac{-t_0^2}{(t_0+1)(s_0t_0+s_0+t_0)}<0,\]
so $u_o$ must be on the left of $u_{[r_0+1]}$, $u_{[s_0+1]}$ and $u_{[t_0+1]}$. Suppose the edgepath of the $[r_0, \!r_1,\ldots,r_{m-1},\! r_{m}]$-tangle ends on the edge $\langle\frac{1}{r_0+k+1}\rangle\text{---}\langle\frac{1}{r_0+k}\rangle$, where $0\leq k\leq -r_0-2$, then $u_0$ must be the $\textit{u}$-coordinate of ending points of the admissible edgepath system $\Gamma$ below:
\[
\begin{split}
\gamma_{1}:&\!\leftharpoondown\![r_0, r_1,\ldots, r_{m-1},\! -1]\!\leftharpoondown\![r_0, r_1,\ldots, r_{m-1}, \!-2] \!\leftharpoondown\!\cdots \!\leftharpoondown\![r_0, \!r_1,\ldots,r_{m-1},\! r_{m}] \\
           &\!\leftharpoondown\![r_0,\! r_1,\!\ldots,\!r_{m-2},\! -1\!]\!\leftharpoondown\![r_0, r_1,\!\ldots,\! r_{m-2},\!-2]\!\leftharpoondown\! \cdots \!\leftharpoondown\![r_0, \!r_1,\!\ldots,\!r_{m-2}\!,\! r_{m-1}\!+\!2\!]\\
           &\cdots\cdots \\
           &\leftharpoondown[r_0,-1]\leftharpoondown[r_0,-2],\ldots, \leftharpoondown[r_0,r_1+2]\\
           &(\!\frac{-s_0t_0}{s_0\!+\!t_0}\!-\!1\!-\!r_0\!-\!k\!)[\!r_0\!+\!k\!+\!1\!]\!+\!(\!\frac{s_0t_0}{s_0\!+\!t_0}\!+\!2\!+\!r_0\!+\!k)[r_0\!+\!k\!]\!\leftharpoonup\![r_0\!+\!k]\!\cdots \!\leftharpoondown\![r_0\!+\!2\!].
\end{split}
\]

\[
\begin{split}
\gamma_{2}:&\!\leftharpoondown\![s_0,\! s_1,\!\ldots,\! s_{p-1},\! -1]\!\leftharpoondown[s_0,\! s_1,\!\ldots\!, s_{p-1}, \!-2]\! \leftharpoondown\!\cdots \!\leftharpoondown\![s_0, s_1,\ldots,s_{p-1}, s_{p}] \\
           &\!\leftharpoondown\![s_0\!, s_1\!,\!\ldots\!, \!s_{p-2}, \!-1]\! \leftharpoondown[s_0,\! s_1\!,\!\ldots\!,\! s_{p-2}, \!-2]\! \leftharpoondown \!\cdots \!\leftharpoondown\![s_0, s_1,\!\ldots\!,s_{p-2}, s_{p-1}\!+2\!]\\
           &\cdots\cdots \\
           &\frac{s_0}{s_0+t_0}\langle 0\rangle+\frac{t_0}{s_0+t_0}[s_0,-1]\leftharpoondown[s_0,-1]\leftharpoondown[s_0,-2],\ldots, \leftharpoondown[s_0,s_1+2].\\
\end{split}
\]

\[
\begin{split}
\gamma_{3}:&\!\leftharpoondown\![t_0,\! t_1,\!\ldots,\! t_{q-1},\! -1]\!\leftharpoondown[t_0,\! t_1,\!\ldots\!, t_{q-1}, \!-2]\! \leftharpoondown\!\cdots \!\leftharpoondown\![t_0, t_1,\ldots,t_{q-1}, t_{q}] \\
           &\!\leftharpoondown\![t_0\!, t_1\!,\!\ldots\!, \!t_{q-2}, \!-1]\! \leftharpoondown[t_0,\! t_1\!,\!\ldots\!,\! t_{q-2}, \!-2]\! \leftharpoondown \!\cdots \!\leftharpoondown\![t_0, t_1,\!\ldots\!,s_{q-2}, t_{q-1}\!+2\!]\\
           &\cdots\cdots \\
           &\frac{t_0}{s_0+t_0}\langle 0\rangle+\frac{s_0}{s_0+t_0}[t_0,-1]\leftharpoondown[t_0,-1]\leftharpoondown[t_0,-2],\ldots, \leftharpoondown[t_0,t_1+2].\\
\end{split}
\]
The length of the partial edges are calculated via $u_0$ by formula (3.1) from~\cite{IM07}.

The cycle of \textit{r}-value of above edgepath system is $(1, s_0, t_0)$. Any candidate surface associated to this edgepath system is essential by Theorem ~\ref{criteria}(1) or (2).

By formula (4.1), the twist of an essential surface $S_1$ associated to the above edgepath system is
\[
\begin{split}
\tau(S_1)&=\sum_{\gamma_{i}\in \Gamma_{non-const}^{}} \sum_{e_{i,j}\in \gamma_{i}}-2\sigma(e_{i,j})|e_{i,j}|\\
       &=2[(-r_m-1)+(-r_{m-1}-2)+\cdots+(-r_1-2)]\\
       &+\frac{2s_0}{s_0+t_0}+2(k-1)+2(\frac{-s_0t_0}{s_0+t_0}-1-r_0-k)\\
       &+2[(-s_p-1)+(-s_{p-1}-2)+\cdots+(-s_1-2)]+ 2\frac{s_0}{s_0+t_0}\\
       &+2[(-t_q-1)+(-t_{q-1}-2)+\cdots+(-t_1-2)]+ 2\frac{t_0}{s_0+t_0}\\
       &=\frac{2t_0^2}{s_0+t_0}\!-\!2(r_0+t_0+\!2)+\!8-\!4(m+p+q)-\!2(\sum_{i=1}^{m}r_i + \!\sum_{j=1}^{p}s_j +\!\sum_{k=1}^{q}t_k).
\end{split}
\]

So the boundary slope of $S_1$ is
\[
\begin{split}bs(S_1)&=\tau(S_1)-\tau(S_0)\\
&=\frac{2t_0^{2}}{s_0+t_0\!}\!-\!2(r_0+\!t_0+\!2)+6-\!2(m+p+q)-\!2(\!\sum_{even}^{m}r_i\! +\! \sum_{even}^{p}s_j\!+\! \sum_{even}^{q}t_k\!).
\end{split}
\]
By Lemma 4.3 (1), we have
\[
\begin{split}
 &-\frac{\chi(S_1)}{\sharp S_1}=\sum_{\gamma_{i}\in \Gamma_{non-const}}|\gamma_{i}|+N_{const}-N+(N-2-\sum_{\gamma_{i}\in \Gamma_{const}}\frac{1}{q_{i}})\frac{1}{1-u_{0}}\\
 &=(-r_m -1)+(-r_{m-1}-2)+\cdots+(-r_1-2)+(k-1)+(\frac{-s_0t_0}{s_0+t_0}-r_0-k-1)\\
 &+(-s_p -1)+(-s_{p-1}-2)+\cdots+(-s_1-2)+\frac{s_0}{s_0+t_0}\\
 &+(-t_q -1)+(-t_{q-1}-2)+\cdots+(-t_1-2)+\frac{t_0}{s_0+t_0}-3+\frac{s_0t_0+s_0+t_0}{s_0+t_0}\\
 &=-r_0-2(m+p+q)-(\sum_{i=1}^{m}r_i + \!\sum_{j=1}^{p}s_j +\!\sum_{k=1}^{q}t_k).
\end{split}
\]

So far we have proved the case (1) of Theorem 2.5.

{\bf (2)} When $A$ and $C< 0$ \ and \ $\Delta\geq0$, we choose the following admissible edgepath system.
\[
\begin{split}
\beta_{1}:&\!\leftharpoondown\![r_0, r_1,\ldots, r_{m-1},\! -1]\!\leftharpoondown\![r_0, r_1,\ldots, r_{m-1}, \!-2] \!\leftharpoondown\!\cdots \!\leftharpoondown\![r_0, \!r_1,\ldots,r_{m-1},\! r_{m}] \\
           &\!\leftharpoondown\![r_0,\! r_1,\!\ldots,\! r_{m-2},\! -1] \!\leftharpoondown\![r_0, r_1,\!\ldots,\! r_{m-2},\! -2] \!\leftharpoondown\! \cdots \!\leftharpoondown\![r_0, r_1,\!\ldots,r_{m-2},\! r_{m-1}\!+\!2]\\
           &\cdots\cdots \\
           &\langle 0\rangle\leftharpoonup[r_0,-1]\leftharpoondown[r_0,-2],\ldots, \leftharpoondown[r_0,r_1+2].\\
\end{split}
\]

\[
\begin{split}
\beta_{2}:&\!\leftharpoondown\![s_0,\! s_1,\!\ldots,\! s_{p-1},\! -1]\!\leftharpoondown[s_0,\! s_1,\!\ldots\!, s_{p-1}, \!-2]\! \leftharpoondown\!\cdots \!\leftharpoondown\![s_0, s_1,\ldots,s_{p-1}, s_{p}] \\
           &\!\leftharpoondown\![s_0\!, s_1\!,\!\ldots\!, \!s_{p-2}, \!-1]\! \leftharpoondown[s_0,\! s_1\!,\!\ldots\!,\! s_{p-2}, \!-2]\! \leftharpoondown \!\cdots \!\leftharpoondown\![s_0, s_1,\!\ldots\!,s_{p-2}, s_{p-1}\!+2\!]\\
           &\cdots\cdots \\
           &\langle 0\rangle\leftharpoondown[s_0,-1]\leftharpoondown[s_0,-2],\ldots, \leftharpoondown[s_0,s_1+2].\\
\end{split}
\]

\[
\begin{split}
\beta_{3}:&\!\leftharpoondown\![t_0,\! t_1,\!\ldots,\! t_{q-1},\! -1]\!\leftharpoondown[t_0,\! t_1,\!\ldots\!, t_{q-1}, \!-2]\! \leftharpoondown\!\cdots \!\leftharpoondown\![t_0, t_1,\ldots,t_{q-1}, t_{q}] \\
           &\!\leftharpoondown\![t_0\!, t_1\!,\!\ldots\!, \!t_{q-2}, \!-1]\! \leftharpoondown[t_0,\! t_1\!,\!\ldots\!,\! t_{q-2}, \!-2]\! \leftharpoondown \!\cdots \!\leftharpoondown\![t_0, t_1,\!\ldots\!,s_{q-2}, t_{q-1}\!+2\!]\\
           &\cdots\cdots \\
           &\langle 0\rangle\leftharpoondown[t_0,-1]\leftharpoondown[t_0,-2],\ldots, \leftharpoondown[t_0,t_1+2].\\
\end{split}
\]

The cycle of \textit{r}-value of above edgepath system is $(-r_0-2, s_0, t_0)$. Any candidate surface associated to this edgepath system is essential by Theorem ~\ref{criteria}(1) or (3).

The twist of an essential surface $S_2$ associated to the above edgepath system is
\[
\begin{split}
\tau(S_2)=&2[(-r_m-1)+(-r_{m-1}-2)+\cdots+(-r_1-2)] -1\\
       &+2[(-s_p-1)+(-s_{p-1}-2)+\cdots+(-s_1-2)]+ 1\\
       &+2[(-t_q-1)+(-t_{q-1}-2)+\cdots+(-t_1-2)]+ 1\\
       &=8-\!4(m+p+q)-\!2(\sum_{i=1}^{m}r_i + \!\sum_{j=1}^{p}s_j +\!\sum_{k=1}^{q}t_k).
\end{split}
\]
The boundary slope of $S_2$ is
\[bs(S_2)=\tau(S_2)-\tau(S_0)=6-2(m+p+q)-2(\sum_{even}^{m}r_i + \sum_{even}^{p}s_j+ \sum_{even}^{q}t_k).\]
By lemma 4.3 (2), we have
\[
\begin{split}
-\frac{\chi(S_2)}{\sharp S_2}&=(-r_m -1)+(-r_{m-1}-2)+\cdots+(-r_1-2)+1\\
 &+(-s_p -1)+(-s_{p-1}-2)+\cdots+(-s_1-2)+1\\
 &+(-t_q -1)+(-t_{q-1}-2)+\cdots+(-t_1-2)+1-2\\
 &=4-2(m+p+q)-(\sum_{i=1}^{m}r_i + \!\sum_{j=1}^{p}s_j +\!\sum_{k=1}^{q}t_k).
\end{split}
\]

\end{proof}

\section*{Acknowledgments}
The authors would like to thank Ying Zhou for her Matlab program to check some cases of Theorem 2.4. Nathan Dunfield's Python program (available at http://dunfield.info/montesinos) is applied to check Theorem 2.5 for some cases. Special thanks go to Miaowang Li for her constant encouragement and numerous helpful suggestions.

\bibliographystyle{amsplain}

\end{document}